\documentclass[12pt]{amsart}
\usepackage{amsmath,amscd,amssymb,amsfonts,graphics}
\newcommand{\h}{\hbox}
\newcommand{\q}{\quad}
\newcommand{\nin}{\noindent}
\newcommand{\bs}{\par\bigskip}
\newcommand{\ms}{\par\medskip}
\newcommand{\sk}{\par\smallskip}
\newcommand{\bsn}{\par\bigskip\noindent}
\newcommand{\msn}{\par\medskip\noindent}
\newcommand{\skn}{\par\smallskip\noindent}
\newcommand{\ssb}{\raise.15ex\h{${\scriptscriptstyle\bullet}$}}
\newcommand{\ssc}{\,\raise.15ex\h{${\scriptstyle\circ}$}\,}

\newcommand{\mopl}{\hbox{$\bigoplus$}}
\newcommand{\C}{{\mathbb C}}
\newcommand{\DD}{{\mathbb D}}
\newcommand{\N}{{\mathbb N}}
\newcommand{\PP}{{\mathbb P}}
\newcommand{\Q}{{\mathbb Q}}
\newcommand{\R}{{\mathbf R}}
\newcommand{\Z}{{\mathbb Z}}
\newcommand{\A}{{\mathbb A}}
\newcommand{\D}{{\mathcal D}}
\newcommand{\E}{{\mathcal E}}
\newcommand{\F}{{\mathcal F}}
\newcommand{\Hc}{{\mathcal H}}
\newcommand{\I}{{\mathcal I}}
\newcommand{\Lc}{{\mathcal L}}

\newcommand{\OO}{{\mathcal O}}
\newcommand{\pp}{{}^{\mathbf p}}
\newcommand{\pH}{{}^{\mathbf p}{\mathcal H}}

\newcommand{\Zt}{\widetilde{Z}}
\newcommand{\Gr}{{\rm Gr}}
\newcommand{\bl}{\bigl}
\newcommand{\br}{\bigr}
\newcommand{\into}{\hookrightarrow}
\newcommand{\simto}{\buildrel\sim\over\longrightarrow}
\newcommand{\onto}{\mathop{\rlap{$\to$}\hskip2pt\h{$\to$}}}
\newcommand{\indlim}{\rlap{\raise-5.5pt\h{$\,\to$}}{\rm lim}}
\newcommand{\mtim}{\,{\times}\,}
\newcommand{\ges}{\geqslant}
\newcommand{\les}{\leqslant}
\newcommand{\1}{\hskip1pt}
\begin{document}
\title[Dependence of Lyubeznik numbers]{Dependence of Lyubeznik numbers of cones of projective schemes on projective embeddings}
\author[T.~Reichelt]{Thomas Reichelt}
\address{T.~Reichelt : Mathematisches Institut, Universit\"at Heidelberg,
Im Neuenheimer Feld 205, 69120 Heidelberg, Germany}
\email{treichelt@mathi.uni-heidelberg.de}
\author[M.~Saito]{Morihiko Saito}
\address{M. Saito : RIMS Kyoto University, Kyoto 606-8502 Japan}
\email{msaito@kurims.kyoto-u.ac.jp}
\author[U.~Walther]{Uli Walther}
\address{U.~Walther : Purdue University, Dept.\ of Mathematics, 150 N.\ University St., West Lafayette, IN 47907, USA}
\email{walther@math.purdue.edu}
\begin{abstract} We construct complex projective schemes with Lyubeznik numbers of their cones depending on the choices of projective embeddings. This answers a question of G.~Lyubeznik in the characteristic 0 case. It contrasts with a theorem of W.~Zhang in the positive characteristic case where the Frobenius endomorphism is used. Reducibility of schemes is essential in our argument. B.~Wang recently constructed examples of irreducible projective schemes (which are not normal) from our examples of reducible ones. So the question is still open in the normal singular case.
\end{abstract}
\maketitle
\bs
\centerline{\bf Introduction}
\bsn
For a Noetherian local ring containing a field, the Lyubeznik numbers are defined by using a surjection from a {\it regular\1} local ring together with the {\it local cohomology,} see \cite{L1}. Here Lyubeznik treated mainly the characteristic 0 case inspired by some finiteness theorem from \cite{HS} in positive characteristic. (If the ring contains only a field which is not canonically isomorphic to the residue field, we need some completion argument together with the Cohen structure theorem.) The theory has been extended to the mixed characteristic case since then, see \cite{NW}, \cite{HNPW}. For a motivic aspect of local cohomology and Lyubeznik numbers, see \cite{Ga}. (For the case of affine cones over the base field $\C$, see (1) below.)
\sk
These numbers are interesting invariants of the local ring with many applications, see for instance \cite{NWZ}. In characteristic 0, Lyubeznik and many other people used $D$-modules to study these numbers (while the Frobenius endomorphism is the main tool in positive characteristic). However, they do not seem to have employed systematically the theory of {\it regular holonomic\1} $D$-modules, especially, the {\it Riemann-Hilbert correspondence\1} in the base field $\C$ case. Here we do not need the equivalence of categories, but the {\it compatibility\1} of the corresponding functors is essential, see \cite{Gr2}, \cite{Bo}, etc.
\sk
Using this compatibility, it is possible to define the Lyubeznik numbers ``topologically" (see Proposition~1 below), and apply the ``derived algebraic topology". (This was partly done in \cite{GS} for the {\it isolated\1} singularity case, for instance, the affine cone of a {\it smooth\1} projective variety.) Here ``derived" means that it involves the theory of {\it derived categories of bounded $\Q$-complexes with constructible cohomology sheaves\1} $D^b_c(X,\Q)$ (see for instance \cite{BBD}, \cite{Di}, \cite{KK}, etc.), which was initiated by the Grothendieck school in order to solve the Weil conjecture in the $\ell$-adic case. This gives a powerful refinement of the classical algebraic topology, for instance, a generalization of the classical Thom-Gysin sequence, which is the key to the proof of the main theorem. Complexes with constructible cohomology sheaves are easier to handle than complexes of quasi-coherent sheaves, and we can get satisfactory computations of the invariants.
\sk
In this paper we consider only the case of affine cones of projective schemes.
Let $X$ be a projective scheme over $\C$ with $\Lc$ a very ample line bundle. Let $C$ be the cone of $X$ associated with $\Lc$. Let $x_1,\dots,x_n$ be projective coordinates of $Y:=\PP^{n-1}_{\C}$ containing $X$ so that $\OO_{\PP^{n-1}}(1)|_X=\Lc$. Let $I$ be the ideal of $R:=\C[x_1,\dots,x_n]$ defining the cone $C\subset\A^n_{\C}$.
The {\it Lyubeznik numbers} $\lambda_{k,j\1}(C)$ are defined by
$$\lambda_{k,j\1}(C):=\dim_{\C}{\rm Ext}_R^k(\C,H^{n-j}_IR)\q\q(k,j\in\N),
\leqno(1)$$
see \cite{L1,L2,NWZ}. Here the $H^{n-j}_IR$ are the local cohomology groups, and vanish for $j>\dim C$, see Remark~(ii) after (1.1) below. The higher extension groups ${\rm Ext}_R^k(\C,*)$ can be calculated by the Koszul complex for the multiplications of the $x_i$ ($i\in[1,n]$) which gives a free resolution of $\C$ over $R$. This holds also for the higher torsion groups ${\rm Tor}^R_{n-k}(\C,*)$. Setting $V:={\rm Spec}\,R=\A^n_{\C}$, we then get the isomorphisms
$${\rm Ext}_R^k(\C,H^{n-j}_IR)={\rm Tor}^R_{n-k}(\C,H^{n-j}_IR)=
H^{k-n}\1{\bf L}\1i^{\1*}_{0,V}(\Hc^{n-j}_C\OO_{V}).
\leqno(2)$$
Here $i_{A,B}:A\into B$ denotes an inclusion of a subset $A\subset B$ in general, ${\bf L}\1i^{\1*}_{0,V}$ means the derived pull-back functor for $\OO$-modules endowed with an integrable connection (that is, for left $D$-modules), and the $\Hc^{n-j}_C$ are the algebraic local cohomology functors for the closed subscheme $C\subset V$. Note that $\OO_{V}$ is algebraic so that $\Gamma(V,\OO_{V})=R$.
\sk
Using the Riemann-Hilbert correspondence together with the scalar extension by $\Q\into\C$, we then get the following (see (1.1) below).
\msn
{\bf Proposition~1.} {\it In the above notation, we have the equalities}
$$\lambda_{k,j\1}(C)=\dim_{\Q}H^ki_{0,C}^{\1!}(\pH^{-j}\DD\Q_C)\q\q(k,j\in\N).$$
\ms
(See also \cite{GS}.) Here $\Q_{C^{\rm an}}$ and its dual $\DD\Q_{C^{\rm an}}$ (see \cite{Ve1}) are respectively denoted by $\Q_C$ and $\DD\Q_C$ to simplify the notation (where $C^{\rm an}$ is the analytic space associated with $C$), and similarly for $\Q_X$, $\DD\Q_X$. The $\pH^j$ are the cohomology functors associated with the truncations $\pp\tau_{\les k}$ constructed in \cite{BBD} (see also \cite{Di,KS}). Note that the usual cohomology functors $\Hc^j$ for bounded complexes of $D$-modules having regular holonomic cohomology sheaves correspond to the functors $\pH^j$ by the Riemann-Hilbert correspondence.
\sk
Setting
$$\F_j:=\pH^{-j}\DD\Q_C,\q\F'_j:=\F_j|_{C'}\q\,\h{with}\q\, C':=C\setminus\{0\}\q(j\in\N),
\leqno(3)$$
we have the following (see (1.2) below).
\msn
{\bf Proposition~2.} {\it For $k\ges 2$, there are isomorphisms}
$$H^ki_{0,C}^{\1!}\1\F_j=H^{k-1}i_{0,C}^*\1\R(i_{C',C})_*\1\F'_j=H^{k-1}(C',\F'_j).$$
\ms
Combined with Propositions~1, this implies the following.
\msn
{\bf Corollary~1.} {\it We have}
$$\lambda_{k,j\1}(C)=\dim_{\Q}H^{k-1}(C',\F'_j)\q\q(k\ges 2).$$
\ms
For $k\in\Z$, $j\in\N$, set
$$\aligned H_{(j)}^k(X)&:=H^k(X,\pH^{-j}\DD\Q_X),\\
H_{(j)}^k(X)^{\Lc}&:={\rm Ker}\bl(c_1(\Lc):H_{(j)}^k(X)\to H_{(j)}^{k+2}(X)(1)\br),\\ H_{(j)}^k(X)_{\Lc}&:={\rm Coker}\bl(c_1(\Lc):H_{(j)}^{k-2}(X)(-1)\to H_{(j)}^k(X)\br),\endaligned
\leqno(4)$$
where $(m)$ denotes a Tate twist for $m\in\Z$, see \cite{De1}.
(Note that a subquotient of $H_{(j)}^k(X)$ is identified with $\Gr_G^jH_{j-k}(X)$ by a spectral sequence, where $\{G^j\}$ is a decreasing filtration on $H_{j-k}(X)=H^{k-j}(X,\DD\Q_X)$ induced by $\pp\tau_{\les -j}$ on $\DD\Q_X$, see \cite{BBD,Ve2}.)
\sk
Using a {\it generalized Thom-Gysin sequence,} we get the following (see (1.3--4) below).
\msn
{\bf Proposition~3.} {\it There are short exact sequences}
$$0\to H_{(j-1)}^k(X)_{\Lc}(1)\to H^{k-1}(C',\F'_j)\to H_{(j-1)}^{k-1}(X)^{\Lc}\to 0\q(k\in\Z).$$
\sk
Combining this with Corollary~1, we get the following.
\msn
{\bf Corollary~2.} {\it The Lyubeznik numbers $\lambda_{k,j\1}(C)$ of the cone $C$ of a projective scheme $X$ depend on the choice of a very ample line bundle $\Lc$ if the following condition holds\,$:$
\skn
\rlap{\rm (5)}\hskip1.5cm\hangindent=1.5cm
$\dim H_{(j-1)}^k(X)_{\Lc}+\dim H_{(j-1)}^{k-1}(X)^{\Lc}\,$ depends on $\,\Lc\,$ for some $k\ges 2$, $j\ges 1$.}
\ms
(Here it is not very clear whether we can replace $H_{(j-1)}^k(X)_{\Lc}$, $H_{(j-1)}^{k-1}(X)^{\Lc}$ with the graded pieces of the weight filtration $W$ unless we consider the corresponding refinement of the Lyubeznik numbers. There might be a cancellation among the dependences of dimensions for various weights. It seems quite non-trivial to apply a semi-continuity argument even after replacing the $\Q$-coefficients with $\C$, since the rank of a morphism may become strictly higher by taking a linear combination.)
The converse of Corollary~2 holds in certain cases, see Corollary~1.7 below. This implies the independence of Lyubeznik numbers under projective embeddings in the $\Q$-homology manifold case (generalizing \cite{Swl} in the non-singular case), see Corollary~1.8 below.
Using Corollary~2, we can prove the following.
\msn
{\bf Theorem~1.} {\it For any field $K$ of characteristic $0$, there are projective schemes over $K$ such that their irreducible components are smooth and the Lyubeznik numbers $\lambda_{k,j\1}(C)$ of their cones $C$ depend on the choices of projective embeddings for some $k\ges 2$. Here $j$ coincides with the dimension of the lowest dimensional irreducible component of $C$, and $X$ can be equidimensional.}
\ms
This answers a question of G.~Lyubeznik \cite{L2} in the characteristic 0 case (see \cite{Swl} for the $X$ non-singular case where the answer is different). Theorem~1 was rather unexpected, since the situation is entirely different in the positive characteristic case where the Frobenius endomorphism can be used, see \cite{Z}, \cite{NWZ}.
\sk
The proof of Theorem~1 can be reduced to the case $K=\C$ provided that examples are defined over $\Q\subset\C$, since the local cohomology is compatible with the base change under an extension of base field. (The latter property follows from an expression of the derived local cohomology functor using a \v Cech complex consisting of localizations, see for instance \cite[Thm.\ A1.3]{Ei}.) Thus the proof is reduced by Corollary~2 to finding complex projective schemes $X$ satisfying condition~(5) and defined over $\Q\subset\C$.
\sk
It is rather easy to construct schemes satisfying the above conditions if the condition $k\ges 2$ in (5) is omitted (where $k$ may be negative), see (2.1) below. In order to satisfy this condition, we need some more construction, where the argument is easier in the non-equidimensional case (see (2.2) below), and we have to use Hodge theory in the equidimensional case (see (2.3) below). In these arguments, {\it reducibility} of schemes is essential. Botong Wang \cite{Wa} recently succeeded in constructing examples of irreducible projective schemes (which are not normal) based on our examples of reducible ones. So the question is still open in the normal singular case.
\sk
The first named author was supported by a DFG Emmy-Noether-Fellowship (RE 3567/1-1) and acknowledges partial support by the project SISYPH: ANR-13-IS01-0001-01/02 and DFG grant HE 2287/4-1 \& SE 1114/5-1. He would like to thank Duco van Straten for a stimulating discussion.
The second named author is partially supported by Kakenhi 15K04816.
The third named author is supported in part by NFS grant DMS-1401392 and by Simons Foundation Collaboration Grant for Mathematicians \#580839. He thanks Nick Switala for stimulating conversations. The authors thank the referees for their useful comments to improve the paper.
\sk
In Section~1 we review generalized Thom-Gysin sequences, study the Lyubeznik numbers in the $X$ $\Q$-homology manifold case, and prove Propositions~1--3.
In Section~2 we prove Theorem~1 by constructing desired examples.
\bs\bs
\vbox{\centerline{\bf 1. Preliminaries}
\msn
In this section we review generalized Thom-Gysin sequences, study the Lyubeznik numbers in the $X$ $\Q$-homology manifold case, and prove Propositions~1--3.}
\msn
{\bf 1.1~Proof of Proposition~1.}
By the Riemann-Hilbert correspondence for algebraic $D$-modules using the de Rham functor DR (see for instance \cite{Bo}), the derived pull-back functor ${\bf L}\1i^{\1*}_{0,V}[-n]$ (explained before Proposition~1) corresponds to $i^{\1!}_{0,V}$, that is,
$${\rm DR}\ssc{\bf L}\1i^{\1*}_{0,V}[-n]=i^{\1!}_{0,V}\ssc{\rm DR},
\leqno(1.1.1)$$
see also Remark~(iii) below (and \cite[Remark after Corollary~2.24]{Sa2}). We have moreover
$${\rm DR}(\OO_V[n])=\DD\C_V\,(=\C_V[2n]),
\leqno(1.1.2)$$
since ${\rm DR}(\OO_V)=\C_V[n]$. Here $\DD$ denotes the dual functor, see \cite{Ve1}. Note that the functor ${\bf L}\1i^{\1*}_{0,V}[-n]$ corresponds to $i^{\1*}_{0,V}$ under the {\it contravariant\1} functor ${\rm Sol}=\DD\ssc{\rm DR}$, see Remark~(i) below and also \cite{Ka}, \cite{KK}, \cite{Me}, \cite[Remark~2.4.15 (3)]{Sa1}. (The equivalence of categories itself is not really needed here.)
\sk
The derived local cohomology functor $\R\Gamma_C$ corresponds to $(i_{C,V})_*i_{C,V}^{\1!}$ (see \cite{Bo}), and we have
$$i_{C,V}^{\1!}\DD\Q_V=\DD\1i_{C,V}^{\1*}\Q_V=\DD\Q_C.
\leqno(1.1.3)$$
So Proposition~1 follows.
\msn
{\bf Remarks.} (i) Let $X$ be a complex manifold (or a smooth complex algebraic variety) of dimension $n$. For a bounded complex $M^{\ssb}$ of left $\D_X$-modules having regular holonomic cohomology sheaves, the de Rham and solution functors can be defined by
$$\aligned{\rm DR}(M^{\ssb})&:=\R\Hc om_{\D_X}(\OO_X,M^{\ssb})[n],\\{\rm Sol}(M^{\ssb})&:=\R\Hc om_{\D_X}(M^{\ssb},\OO_X)[n].\endaligned
\leqno(1.1.4)$$
(If $X$ is a smooth complex algebraic variety, $X$ and $M^{\ssb}$ on the right-hand side are respectively replaced by $X^{\rm an}$ and $M^{\rm an\,\ssb}$.)
Taking the composition, we can get a perfect pairing
$${\rm DR}(M^{\ssb})\otimes_{\C}{\rm Sol}(M^{\ssb})\to\R\Hc om_{\D_X}(\OO_X,\OO_X)[2n]=\C_X[2n]=\DD\C_X,
\leqno(1.1.5)$$
that is, the corresponding morphism
$${\rm DR}(M^{\ssb})\to\DD\bl({\rm Sol}(M^{\ssb})\br):=\R\Hc om_{\C_X}\bl({\rm Sol}(M^{\ssb}),\DD\C_X\br)
\leqno(1.1.6)$$
is an isomorphism. It is also known that ${\rm DR}$ commutes with $\DD$. (This follows, for instance, from \cite[Proposition~2.4.12]{Sa1}.)
\ms
(ii) In the notation of the introduction, we have by \cite[Theorem 3.8]{Ha1} (see also \cite{BS,Iy})
$$H^{n-j}_IR=0\q\h{for}\q n-j<{\rm codim}\,C.
\leqno(1.1.7)$$
In our case this can be shown by taking a complete intersection $C'$ containing $C\subset\C^n$ with $\dim C'=\dim C$, and using the composition of derived local cohomology functors $\R\Gamma_{C'}$ and $\R\Gamma_C$. Indeed, (1.1.7) is easy in the complete intersection case by the theory of regular sequences (see for instance \cite[Thm.~II.8.21A(c)]{Ha2}), and $\R\Gamma_C$ is the {\it right\1} derived functor of a {\it left exact\1} functor $\Gamma_C$.
\sk
Note also that, by the Riemann-Hilbert correspondence and (1.1.2), the above vanishing is equivalent to
$$\pH^{-j}\1\R\Gamma_C\1\DD\Q_V=0\,\,\,\,\h{(that is,}\,\,\,\pH^j\Q_C=0\,)\q\h{for}\q j>\dim C,
\leqno(1.1.8)$$
using (1.1.3) and a remark before it. The assertion for $\Q_C$ follows easily from the definition of the $t$-structure in \cite{BBD}: we have $K^{\ssb}\in D^b_c(X)^{\les k}$, that is, $\pH^jK^{\ssb}=0$ ($j>k$), if and only if $\dim{\rm supp}\,\Hc^jK^{\ssb}\les k\1{-}\1j$ ($j\in\Z$).
\ms
(iii) Let $X$ be a closed submanifold of a complex manifold $Y$. For a regular holonomic right $\D_Y$-module $M$, there are natural inclusions (as $\OO_Y$-modules)
$$\E xt^j_{\OO_Y}(\OO_X,M)\into\Hc^j_{[X]}M\q(j\in\N),
\leqno(1.1.9)$$
inducing isomorphisms of right $\D_Y$-modules
$$\E xt^j_{\OO_Y}(\OO_X,M)\otimes_{\D_X}\D_{X\to Y}\simto\Hc^j_{[X]}M\q(j\in\N).
\leqno(1.1.10)$$
Here $\D_{X\into Y}:=\OO_X\otimes_{\OO_Y}\D_Y$, and the $\Hc^j_{[X]}M$ are the algebraic local cohomology sheaves defined by
$$\Hc^j_{[X]}M:=\rlap{\raise-10pt\h{$\,\,\,\scriptstyle k$}}\indlim\,\E xt^{\1 j}_{\OO_Y}(\OO_Y/\I_X^k,M)\q(j\in\N),
\leqno(1.1.11)$$
with $\I_X\subset\OO_Y$ the ideal sheaf of $X\subset Y$, see \cite{KK}. Note that the sources of (1.1.9) and (1.1.10) are respectively the cohomological pull-back of $M$ as right $\D$-module by the inclusion $i_{X,Y}:X\into Y$ and its direct image as right $\D$-module by $i_{X,Y}$. (Using a spectral sequence, the proof of (1.1.10) can be reduced to the codimension 1 case.)
\sk
Set $r:={\rm codim}_YX$. The formula corresponding to (1.1.10) for a regular holonomic {\it left} $\D$-module $M$ is as follows (see, for instance, the proof of \cite[Corollary 5.4.6]{KK}):
$$\D_{Y\leftarrow X}\otimes_{\D_X}{\mathcal T}\!or_{r-j}^{\OO_Y}(\OO_X,M)=\Hc_{[X]}^jM\q(j\in\N).
\leqno(1.1.12)$$
\msn
{\bf 1.2~Proof of Proposition~2.}
The last isomorphism in Proposition~2 holds, since $0$ is the vertex of the cone $C$.
The first isomorphism follows from the long exact sequence associated with the distinguished triangle
$$i_{0,C}^{\1!}\1\F_j\to i_{0,C}^{\1*}\1\F_j\to i_{0,C}^{\1*}\1\R(i_{C',C})_*\1 i_{C',C}^*\1\F_j\buildrel{+1}\over\to,
\leqno(1.2.1)$$
since $H^ki_{0,C}^{\1*}\1\F_j=0$ for $k>0$, see \cite{BBD} (and also \cite[Remark after Corollary~2.24]{Sa2}). The last triangle can be obtained by applying the functor $i_{0,C}^{\1*}$ to
the triangle
$$(i_{0,C})_*i_{0,C}^{\1!}\to id\to\R(i_{C',C})_*\1 i_{C',C}^*\buildrel{+1}\over\to.$$
This finishes the proof of Proposition~2.
\msn
{\bf 1.3.~Generalized Thom-Gysin sequences} (see also \cite{Ko, Swz}). Let $\pi:E\to X$ be a vector bundle of rank $r$ on a complex analytic space $X$ which is assumed connected. Set $E':=E\setminus X$, where $X$ is identified with the zero section of $E$. There are natural morphisms
$$i_X:X\into E,\q j_{E'}:E'\into E,\q\pi':=\pi|_{E'}:E'\to X.$$
For $\F^{\ssb}\in D^b_c(X,\Q)$, we have the distinguished triangle
$$\F^{\ssb}\buildrel{\xi}\over\to\F^{\ssb}(r)[2r]\to\R\pi'_*\pi'{}^{\1!}\F^{\ssb}\buildrel{+1}\over\to,
\leqno(1.3.1)$$
inducing a long exact sequence called a {\it generalized Thom-Gysin sequence}\,:
$$\to H^k(X,\F^{\ssb})\buildrel{\xi}\over\to H^{k+2r}(X,\F^{\ssb})(r)\to H^k(E',\pi'{}^{\1!}\F^{\ssb})\to H^{k+1}(X,\F^{\ssb})\to.
\leqno(1.3.2)$$
\sk
Indeed, the triangle (1.3.1) is identified with the distinguished triangle
$$i_X^{\1!}\pi^{\1!}\F^{\ssb}\to\R\pi_*\pi^{\1!}\F^{\ssb}\to\R\pi'_*\pi'{}^{\1!}\F^{\ssb}\buildrel{+1}\over\to,
\leqno(1.3.3)$$
since $\pi^{\1!}\F^{\ssb}=\pi^{-1}\F^{\ssb}(r)[2r]$. The last triangle is obtained by applying $\R\pi_*$ to
$$(i_X)_*i_X^{\1!}\pi^{\1!}\F^{\ssb}\to\pi^{\1!}\F^{\ssb}\to\R(j_{E'})_*j_{E'}^*\pi^{\1!}\F^{\ssb}\buildrel{+1}\over\to.
\leqno(1.3.4)$$
The morphism $\xi$ in (1.3.1) is induced by the Euler class of $E$ via the adjunction isomorphism for $a_X^*$ and $(a_X)_*:$
$${\rm Hom}(\Q_X,\Q_X(r)[2r])={\rm Hom}\bl(\Q,\R\Gamma(\Q_X(r)[2r])\br)=H^{2r}(X,\Q)(r).
\leqno(1.3.5)$$
Here ${\rm Hom}$ denotes the group of morphisms in the derived categories, and the {\it Euler class} of $E$, denoted by $e$, is the image of 1 by the following morphism induced by $\xi$ for $\F^{\ssb}=\Q_X$\.:
$$\Q=H^0(X,\Q)\to H^{2r}(X,\Q)(r),
\leqno(1.3.6)$$
which is identified with an element in the first term of (1.3.5), see also \cite[Ex. III.7]{KS}.
\sk
Note that the Euler class $e$ of $E$ induces conversely the morphism $\xi$ in (1.3.1) by taking the tensor product of $e$ (as an element in the first term of (1.3.5)) with the identity on $\F^{\ssb}$\,:
$$\F^{\ssb}=\Q_X\otimes_{\Q}\F^{\ssb}\buildrel{e\otimes id}\over\longrightarrow\Q_X(r)[2r]\otimes_{\Q}\F^{\ssb}=\F^{\ssb}(r)[2r].
\leqno(1.3.7)$$
Indeed, $\xi$ coincides with the tensor product of $\xi$ for $\F^{\ssb}=\Q_X$ (that is, the Euler class $e$ of $E$) with the identity on $\F^{\ssb}$ as is seen from the commutative diagram
$$\begin{array}{ccc}
i_X^{\1!}\pi^{-1}\Q_X\otimes_{\Q}\F^{\ssb}&\to&\R\pi_*\pi^{-1}\Q_X\otimes_{\Q}\F^{\ssb}\\ \downarrow&\raise15pt\h{}\raise-10pt\h{}&\downarrow\,\\ i_X^{\1!}\pi^{-1}\F^{\ssb}&\to&\R\pi_*\pi^{-1}\F^{\ssb}\end{array}\leqno(1.3.8)$$
since (1.3.1) comes from (1.3.3).
This diagram is shown by taking $\pi_*$ of the commutative diagram
$$\begin{array}{ccccc}
(\Gamma_XI^{\ssb})\otimes_{\Q}\pi^{-1}\F^{\ssb}&\into&\Gamma_X(I^{\ssb}\otimes_{\Q}\pi^{-1}\F^{\ssb})&\into&I^{\ssb}\otimes_{\Q}\pi^{-1}\F^{\ssb}\\ &&\downarrow&\raise15pt\h{}\raise-10pt\h{}&\downarrow\,\\ &&\Gamma_XJ^{\ssb}&\into&J^{\ssb}\end{array}$$
since $\pi:E\to X$ is locally trivial on $X$. Here $I^{\ssb}$ is a flasque resolution of $\Q_E=\pi^{-1}\Q_X$, and $J^{\ssb}$ is a flasque resolution of $I^{\ssb}\otimes\pi^{-1}\F^{\ssb}$ (and $\Gamma_X$ denotes the local cohomology functor associating the subsheaf of local sections supported in $X\subset E$ as usual).
\msn
{\bf Remark.} The long exact sequence (1.3.2) holds in the category of mixed $\Q$-Hodge structures if $\F^{\ssb}$ underlies a bounded complex of mixed Hodge modules. Indeed, the above construction can be lifted naturally in the category of mixed Hodge modules, see \cite{Sa2}. For the calculation of $\xi$, we can use the faithfulness of the forgetful functor associating the underlying $\Q$-vector space of a mixed Hodge structure.
\msn
{\bf 1.4.~Proof of Proposition~3.} In the notation of (1.3) and the introduction, we have
$$\pi'{}^{\1!}\1\F^{\ssb}=\F'_j{}^{\ssb}[1]\q\h{by setting}\q\F^{\ssb}:=\pp\Hc^{-j+1}\DD\Q_X,
\leqno(1.4.1)$$
where $E'=C'$ and $r=1$. So Proposition~3 follows from (1.3.2).
\msn
{\bf 1.5.~$\Q$-homology manifold case} (see \cite{GS} for the $X$ smooth case). Let $X$ be a projective scheme such that
$$\pH^j\Q_X=0\q(j\ne d),
\leqno(1.5.1)$$
where $d\in\Z_{>0}$. This condition is satisfied for instance if $X^{\rm an}$ is purely $d$-dimensional, and is a $\Q$-homology manifold or analytic-locally a complete intersection. (The proof of the last assertion follows, for instance, by using the Riemann-Hilbert correspondence and the local cohomology sheaves defined as in (1.1.11).)
\sk
In the notation of the introduction, the assumption~(1.5.1) implies that
$${\rm Supp}\,\pH^k\1\DD\Q_C,\,\,{\rm Supp}\,\pH^k\1\R(j_{C'})_*\DD\Q_{C'}\,\subset\,\{0\}\q(k\ne-d-1).
\leqno(1.5.2)$$
We have the distinguished triangle
$$\Q_{\{0\}}\to\DD\Q_C\to\R (j_{C'})_*\DD\Q_{C'}\buildrel{+1}\over\to,
\leqno(1.5.3)$$
which is the dual of the short exact sequence
$$0\to (j_{C'})_!\Q_{C'}\to\Q_C\to\Q_{\{0\}}\to 0.$$
In this section $j_{C',C}$ and $i_{0,C}$ are denoted respectively by $j_{C'}$ and $i_0$ to simplify the notation.
\sk
We have the Leray-type spectral sequences
$$^*\!E_2^{p,q}=H^pi_0^{\1*}\,\pH^q\1\R (j_{C'})_*\DD\Q_{C'}\Longrightarrow H^{p+q}i_0^{\1*}\1\R (j_{C'})_*\DD\Q_{C'},
\leqno(1.5.4)$$
\vskip-6.5mm
$$\q\q^!\!E_2^{p,q}=H^pi_0^{\1!}\,\pH^q\1\R (j_{C'})_*\DD\Q_{C'}\Longrightarrow H^{p+q}i_0^{\1!}\1\R (j_{C'})_*\DD\Q_{C'}=0,
\leqno(1.5.5)$$
where the last vanishing follows from $i_0^{\1!}\1\R (j_{C'})_*=0$ (the latter is the dual of $i_0^{\1*}\1\R (j_{C'})_!=0$).
These can be constructed by using spectral objects in \cite{Ve2} together with an argument similar to \cite[Example 1.4.8]{De1}, or we can use the Riemann-Hilbert correspondence after the scalar extension by $\Q\into\C$.
\sk
By (1.1.8), (1.5.2) together with properties of $i_0^{\1*}$, $i_0^{\1!}$ in \cite{BBD} (see also \cite[Remark after Corollary~2.24]{Sa2}), we get
$${}^*\!E_2^{-p,q}={}^!\!E_2^{p,q}=0\q\h{unless}\,\,\,\,p=0,\,q\ges -d-1\,\,\,\,\h{or}\,\,\,\,q=-d-1,\,\,p\ges 0.
\leqno(1.5.6)$$
The generalized Thom-Gysin sequence (1.3.2) together with an isomorphism similar to the last isomorphism of Proposition~2 implies that
$$H^{-k}i_0^{\1*}\1\R (j_{C'})_*\DD\Q_{C'}=0\q\h{unless}\q k\in[1,2d+2].
\leqno(1.5.7)$$
(This can be shown also by using the link $L_{C,0}$ of $C$ at $0$. It is the intersection of $C$ with a sphere $S^{2n-1}$ around $0\in\C^n$, and is a $(2d{+}1)$-dimensional real analytic space, see \cite{DS}. Its dualizing complex $\DD\Q_{L_{C,0}}$ (see \cite{Ve1}) is isomorphic to the restriction of $\DD\Q_{C'}[-1]$ to $L_{C,0}$.)
\sk
The spectral sequence (1.5.4) degenerates at $E_2$ by (1.5.6), and it follows from (1.1.8), (1.5.2), (1.5.7) that
$$\pH^{-j}\1\R (j_{C'})_*\DD\Q_{C'}=0\q\h{unless}\q j\in[1,d+1].
\leqno(1.5.8)$$
\sk
Using (1.5.3) and (1.5.8), we can prove the isomorphisms
$$\pH^k\1\DD\Q_C\simto\pH^k\1\R(j_{C'})_*\DD\Q_{C'}\q(k\ne -1),
\leqno(1.5.9)$$
together with the short exact sequence
$$0\to\pH^{-1}\DD\Q_C\to\pH^{-1}\R(j_{C'})_*\DD\Q_{C'}\to\Q_{\{0\}}\to 0\,.
\leqno(1.5.10)$$
Here the vanishing of $\pH^0\DD\Q_C$ is rather nontrivial. If we have $\pH^0\DD\Q_C\ne 0$, then (1.5.8) and the long exact sequence associated with (1.5.3) imply the surjectivity of the composition
$$\Q_{\{0\}}\to\DD\Q_C\to\pH^0\DD\Q_C.$$
We then get a splitting of the first morphism, but this is a contradiction. So the vanishing of $\pH^0\DD\Q_C$ follows.
\sk
Combined with (1.5.6), the spectral sequence (1.5.5) implies the isomorphisms
$$H^pi_0^{\1!}\,\pH^{-d-1}\1\R (j_{C'})_*\DD\Q_{C'}=\begin{cases}\pH^{p-d-2}\1\R (j_{C'})_*\DD\Q_{C'}&\h{if}\q p\ges 2,\\ \,0&\h{if}\q p=0,1.\end{cases}
\leqno(1.5.11)$$
Here the direct image $(i_0)_*$ is omitted on the left-hand side to simplify the notation.
\sk
By (1.5.9--11), we get the following (see \cite{GS} for the $X$ smooth case).
\msn
{\bf 1.6.~Proposition.} {\it Under the assumption~$(1.5.1)$, we have
$$\lambda_{k,j}(C)=0\q\h{unless}\q j=d+1,\,\,k\in[2,d+1]\q\h{or}\q k=0,\,\,j\in[1,d],
\leqno(1.6.1)$$
and moreover the following relations among the Lyubeznik numbers hold\,$:$
$$\lambda_{k,d+1}(C)=\lambda_{0,d+2-k}(C)+\delta_{k,d+1}\q(k\in[2,d+1]),
\leqno(1.6.2)$$
where $\delta_{k,d+1}=1$ if $\,k=d+1$, and $0$ otherwise.}
\ms
This implies the following.
\msn
{\bf 1.7.~Corollary.} {\it Under the assumption~$(1.5.1)$, the converse of Corollary~$1$ holds.}
\ms
We then get the following generalization of \cite{Swl} in the $X$ non-singular case.
\msn
{\bf 1.8.~Corollary.} {\it The Lyubeznik numbers $\lambda_{k,j}(C)$ of the cone $C$ of a complex projective scheme $X$ are independent of the choice of a projective embedding of $X$, if the associated analytic space $X^{\rm an}$ is a $\Q$-homology manifold.}
\msn
{\it Proof.} By the definition of $\Q$-homology manifold, we have $H^j_{\{x\}}\Q_X=\Q$ if $j=2\dim X$, and $0$ otherwise ($\forall\,x\in X^{\rm an}$). By induction on strata, we see that the composition of the following two canonical morphisms is an isomorphism (see \cite{GM,BBD}):
$$\Q_X[\dim X]\to{\rm IC}_X\Q\to\DD\Q_X(-\dim X)[-\dim X],
\leqno(1.8.1)$$
where the last morphism is the dual of the first. This implies that $\pH^j(\Q_X[\dim X])=0$ ($j\ne0$), and the above two morphisms are both isomorphisms by \cite{GM} or using the simplicity of the intersection complex ${\rm IC}_X\Q$. Moreover the hard Lefschetz theorem holds for the intersection cohomology of the complex projective variety $X$, see \cite{BBD} (and \cite{Sa1}). The dimensions of the kernel and cokernel of the action of $c_1(\Lc)$ on the (intersection) cohomology are then read off from the Betti numbers of $X$ by using the primitive decomposition. So the assertion follows from Corollary~1.7. This finishes the proof of Corollary~1.8.
\msn
{\bf Remark.} The independence of the Lyubeznik numbers of cones also holds if $X$ has only isolated singularities and (1.5.1) is satisfied.
\bs\bs
\vbox{\centerline{\bf 2. Construction of examples}
\bsn
In this section we prove Theorem~1 by constructing desired examples.}
\msn
{\bf 2.1.~Examples with condition~(5) satisfied for ``some $k\in\Z\,$".}
Let $Y$ be a smooth complex projective variety of dimension $d_1\ges 2$ having {\it very ample\1} divisors $D$, $D'$ such that their Chern classes $c_1(D)$, $c_1(D')$ are linearly independent. (The last condition can be satisfied in case the Picard number of $Y$ is at least 2, for instance, if $Y$ is $\PP^1\times\PP^1$ or a one point blow-up of $\PP^2$.)
Consider the line bundle $L$ (as a variety) corresponding to $D$. (In this paper a line bundle usually means an {\it invertible sheaf}.)
We have the section $Y_1$ of $L$ corresponding to $D$. Here we may assume $D$ is effective, reduced and smooth. Then $Y_1$ intersects the zero section $Y_0$ of $L$ {\it transversally\1} along $D$. Set $X_1:=Y_0\cup Y_1\subset L$. This is a divisor with normal crossing, and we get
$$\pH^{d_1}\1\Q_{X_1}=\Q_{X_1}[d_1],\q\pH^{-d_1}\1\DD\Q_{X_1}=\DD\Q_{X_1}[-d_1].
\leqno(2.1.1)$$
Using the dual of the short exact sequence
$$0\to\Q_{X_1}\to\Q_{Y_0}\oplus\Q_{Y_1}\to\Q_D\to 0,$$
which is identified with a distinguished triangle, we get the long exact sequence
$$\aligned&\to H^{j-2}(D)(-1)\buildrel{i_*}\over\to H^j(Y)\oplus H^j(Y)\to H^{j-2d_1}\bl(X_1,\DD\Q_{X_1}(-d_1)\br)\\&\to  H^{j-1}(D)(-1)\to\cdots,\endaligned
\leqno(2.1.2)$$
inducing the isomorphisms
$$\aligned\Gr^W_2H^{2-2d_1}\bl(X_1,\DD\Q_{X_1}(-d_1)\br)&={\rm Coker}\bl(H^0(D)(-1)\buildrel{i_*}\over\into H^2(Y)\oplus H^2(Y)\br),\\H^{-2d_1}\bl(X_1,\DD\Q_{X_1}(-d_1)\br)&=H^0(Y)\oplus H^0(Y).\endaligned
\leqno(2.1.3)$$
Here $i_*$ denotes the Gysin morphism for the inclusion $i:D\into Y$, and $\dim H^0(D)=1$ by the weak Lefschetz theorem.
Since $X_1$ is finite over $Y$, and a finite morphism is {\it ample\1} in the
sense of Grothendieck with relatively ample line bundle {\it trivial,} the ample line bundles corresponding to $D,D'$ on $Y$ give ample line bundles on $X_1$ via the pull-back by the natural morphism $\pi_Y:X_1\to Y$, see \cite[Propositions 4.4.10 and Corollary 6.1.11]{Gr1} and Remark~(i) below.
\sk
By (2.1.3) (together with (2.1.1)), we get a difference in the dimension of the images of
$$c_1(D),\,c_1(D'):H^{-d_1}(X_1,\pH^{-d_1}\DD\Q_{X_1})(-1)\to H^{2-d_1}(X_1,\pH^{-d_1}\DD\Q_{X_1}).
\leqno(2.1.4)$$
Here the action of $\ell\in H^2(X_1)(1)$ on $\F\in D^b_c(X_1,\Q)$ can be defined by taking the tensor product of $\F$ with the morphism $\Q_{X_1}\to\Q_{X_1}(1)[2]$ defined by $\ell$. (This is compatible with the normalization.)
\sk
As a conclusion, condition~(5) is satisfied (except for the condition $k\ges 2$) with
$$j-1=d_1,\q k=2-d_1\,(<2).
\leqno(2.1.5)$$
if we further assume the vanishing of $H^1(Y)$ so that $H^{k-1}_{(j-1)}(X_1)=0$.
\sk
This construction will be used in (2.2--3) below, where $k$ can be ``shifted" by taking a {\it product with an appropriate projective scheme\1} (and using the Segre embedding) so that the condition $k\ges 2$ in (5) will be satisfied .
\msn
{\bf Remarks.} (i) The pull-back to $X_1$ of the very ample line bundle $\Lc_Y:=\OO_Y(D)$ is very ample, since $D$ is assumed smooth.
Indeed, there is a short exact sequence of $\OO_{X_1}$-modules
$$0\to\OO_{X_1}\to\OO_{Y_0}\oplus\OO_{Y_1}\to\OO_D\to 0.$$
Taking the tensor product with $\Lc_1:=\pi_Y^*\Lc_Y$, we can deduce that
$$\Gamma(X_1,\Lc_1)={\rm Ker}\bl(\Gamma(Y,\Lc_Y)\oplus\Gamma(Y,\Lc_Y)\buildrel{\delta}\over\to\Gamma(D,\Lc_Y|_D)\br).$$
Here $\delta(s,s'):=(s-s')|_D$, and $Y_0,Y_1$ are identified with $Y$.
We will denote by $\sigma(s,s')$ the global section of $\Lc_1$ corresponding to $(s,s')\in{\rm Ker}\,\delta$.
Then $\Gamma(X_1,\Lc_1)$ contains $\sigma(s+s'',s)$ for $s,s''\in\Gamma(Y,\Lc_Y)$, where $s''$ is a defining section of $D$. This implies an embedding of $X_1$ into a projective space.
Indeed, its restrictions to $Y_0,Y_1$ are embeddings using the inclusion
$$\Gamma(Y,\Lc_Y)\ni s\mapsto\sigma(s,s)\in\Gamma(X_1,\Lc_1),$$
which corresponds to a {\it projection\1} of projective spaces.
Here the images of $Y_0\setminus D$ and $Y_1\setminus D$ are separated by $s''$ in $\sigma(s+s'',s)$ (combined with the above projection). Moreover, $\sigma(s'',0)/\sigma(s,s)$ gives a local coordinate of $Y_0$ {\it vanishing on\1} $Y_1\subset X_1$, if $s$ does not vanish at a given point (and similarly for $\sigma(0,-s'')/\sigma(s,s)$ with $Y_0$, $Y_1$ exchanged).
\sk
The pull-back of $D'$ is also very ample if $D'-D$ is effective and {\it base-point-free,} where $s''$ is replaced by the product of $s''$ with a section of $\OO_Y(D'-D)$ (which is generated by global sections because of the last assumption).
\ms
(ii) Assume $Y=\PP^1\times\PP^1$ with $D\subset\PP^1\times\PP^1$ the diagonal so that $b_1(D)=0$. Note that $D$ is very ample (since it is linearly equivalent to $E_1+E_2$ with $E_i$ ($i=1,2$) the pull-back of a point of $\PP^1$ by the $i\1$the projection $\PP^1\times\PP^1\to\PP^1$), where we use the Segre embedding $\PP^1\times\PP^1\into\PP^3$.
For $D'$, we can take $a_1E_1+a_2E_2$ for any $a_1,a_2\in\Z_{>0}$ with $a_1\ne a_2$; for instance, $(a_1,a_2)=(2,1)$.
Here we combine the Veronese embedding and the Segre embedding.
\sk
Let $\gamma,\gamma'$ be the dimensions of the cokernels of $c_1(D)$, $c_1(D')$ in (2.1.4) respectively.
Then $d_1=2$, $b_1(Y)=0$, $b_2(Y)=2$, and
$$b_4(X_1)=2,\q b_2(X_1)=3,\q\gamma=2,\q\gamma'=1.
\leqno(2.1.6)$$
\msn
{\bf 2.2.~Non-equidimensional examples with condition~(5) strictly satisfied.} Let $\Zt$ be the blow-up of $\PP^{\1d_2+2}$ along a point $P\in\PP^{\1d_2+2}$, where $d_2\ges 2$. This is identified with a $\PP^1$-bundle over $\PP^{\1d_2+1}$, and we have the projection
$$\rho:\Zt\to\PP^{\1d_2+1},$$
where the target is identified with the set of lines of $\PP^{\1d_2+2}$ passing through $P$. The projection $\rho$ has the zero-section $\Sigma_0$ given by the exceptional divisor of the blow-up. It has another section $\Sigma_{\infty}$ which is disjoint from $\Sigma_0$, and is given by the inverse image of a hyperplane of $\PP^{\1d_2+2}$ not containing the center of the blow-up $P$. Let $W\subset\PP^{\1d_2+1}$ be a linear subspace of codimension 2.
Put
$$Z_1:=\rho^{-1}(W),\q Z_2:=\Sigma_0\sqcup\Sigma_{\infty}\subset\Zt,$$
with
$$\dim Z_1=d_2,\q\dim Z_2=d_2+1.$$
Set
$$X_2:=Z_1\cup Z_2\subset\Zt,\q Z'_1:=Z_1\setminus Z_2.$$
Let $j_{Z'_1}:Z'_1\into Z_1$ be the natural inclusion.
We have the short exact sequence
$$0\to(j_{Z'_1})_!\Q_{Z'_1}\to\Q_{X_2}\to\Q_{Z_2}\to 0.
\leqno(2.2.1)$$
Taking the dual, we get the distinguished triangle
$$\DD\Q_{Z_2}\to\DD\Q_{X_2}\to\R(j_{Z'_1})_*\DD\Q_{Z'_1}\buildrel{+1}\over\to.
\leqno(2.2.2)$$
Note that $Z_1\cap Z_2$ is a divisor on $Z_1$, and $j_{Z'_1}:Z'_1\into Z_1$ is an affine open embedding so that
$$\pH^{-j}\R(j_{Z'_1})_*\DD\Q_{Z'_1}=0\q(j\ne d_2).
\leqno(2.2.3)$$
Since $Z_1,$ $Z_2$ are smooth, we then get
$$\pH^{-j}\DD\Q_{X_2}=\begin{cases}\R(j_{Z'_1})_*\DD\Q_{Z'_1}[-d_2]&\h{if}\,\,\,j=d_2\,,\\ \DD\Q_{Z_2}[-d_2-1]&\h{if}\,\,\,j=d_2+1\,,\\ \,0&\h{otherwise.}\end{cases}
\leqno(2.2.4)$$
This implies that
$$H^k(X_2,\pH^{-d_2}\DD\Q_{X_2})=H^{k-d_2}(Z'_1,\DD\Q_{Z'_1})\cong\begin{cases}\Q&\h{if}\,\,\,k=-d_2,\\ \Q&\h{if}\,\,\,k=d_2-1,\\\,0&\h{otherwise,}\end{cases}
\leqno(2.2.5)$$
where the Tate twist is omitted to simplify the notation. Indeed, $\DD\Q_{Z'_1}=\Q_{Z'_1}(d_2)[2d_2]$ with
$$Z'_1\cong\C^{d_2}\setminus\{0\},$$
since $\Zt\setminus Z_2=\C^{d_2+2}\setminus\{0\}$ which is a $\C^*$-bundle over $\PP^{\1d_2+1}$ by the natural projection, and $Z'_1$ is its restriction over the linear subspace $W\subset\PP^{\1d_2+1}$.
\sk
For $X_1$, $d_1$ as in (2.1) with $H^1(Y)=0$, set
$$X:=X_1\times X_2,\q d:=d_1+d_2+1\,(=\dim X).$$
Then 
$$\pH^{-j}\DD\Q_X=\begin{cases}\DD\Q_{X_1}[-d_1]\boxtimes\R(j_{Z'_1})_*\DD\Q_{Z'_1}[-d_2]&\h{if}\,\,\,j=d-1\,,\\ \DD\Q_{X_1}[-d_1]\boxtimes\DD\Q_{Z_2}[-d_2-1]&\h{if}\,\,\,j=d\,,\\ \,0&\h{otherwise.}\end{cases}
\leqno(2.2.6)$$
\sk
We have very ample line bundles $\Lc_1,\Lc'_1$ on $X_1$ by the pull-backs of $D,D'$ as is explained in (2.1), where $D'-D$ is assumed effective and base-point-free. We choose a very ample line bundle $\Lc_2$ on $X_2$.
These give very ample line bundles $\Lc$, $\Lc'$ on $X$ corresponding to the {\it Segre embedding,} since the very ample line bundles on $X$ are obtained by the {\it tensor product\1} of the pull-backs of the very ample line bundles by the first and second projections from $X=X_1\times X_2$.
Condition~(5) then holds by (2.1.4) and (2.2.4--6) with
$$j-1=d-1,\q k=(2-d_1)+(d_2-1)=d_2-d_1+1,
\leqno(2.2.7)$$
assuming $d_2>d_1$ so that $k\ges 2$. By (2.2.5--6) together with the K\"unneth formula, the last assumption implies that
$$H^k_{(d-1)}(X)\cong\begin{cases}H^{k+d_2-d_1}(X_1,\DD\Q_{X_1})&\h{if}\,\,\,\,|k+d_2|\les d_1,\\ H^{k-d_2+1-d_1}(X_1,\DD\Q_{X_1})&\h{if}\,\,\,\,|k-d_2+1|\les d_1,\\ \,0&\h{otherwise.}\end{cases}
\leqno(2.2.8)$$
So the action of $c_1(\Lc_2)$ vanishes, and the assertion is reduced to the study of the actions of $c_1(\Lc_1)$, $c_1(\Lc'_1)$ in (2.1).
\msn
{\bf Remarks.} (i) Assume $d_2=3$, and $Y,D,D'$ are as in Remark after (2.1), in particular, $d_1=2$. Then we have
$$j=d=6,\q k=2,$$
with
$$H^2_{(5)}(X)\cong H^{-2}(X_1,\DD\Q_{X_1})=H_2(X_1).$$
By Corollary~1 and Proposition~3 together with (2.2.8) and (2.1.6), we see that the Lyubeznik number $\lambda_{2,6}(C)$ is equal to $2$ or $1$ depending on whether the very ample line bundle is induced from $D$ or $D'$.
\sk
(ii) A similar argument holds replacing $Z_1$ with the inverse image of a higher codimensional linear subspace of $\PP^{\1d_2+1}$.
\msn
{\bf 2.3. Equidimensional examples with condition~(5) strictly satisfied.} For an integer $d_2>2$, set
$$B:=\PP^2\times\PP^{\1d_2-2}.$$
Let $\rho':Z'\to B$ be the pull-back of the very ample line bundle on $\PP^{\1d_2-2}$ corresponding to $\OO_{\PP^{\1d_2-2}}(1)$. We have the associated $\PP^1$-bundle
$$\rho:Z\to B.$$
This is a compactification of $Z'$ such that $Z\setminus Z'$ is the section at infinity of $\rho$. Let $Z_2$ be the union of the zero-section and the section at infinity of $\rho$ so that
$$Z_2\,\cong\,\PP^2\mtim\PP^{\1d_2-2}\,\sqcup\,\PP^2\mtim\PP^{\1d_2-2}\,\subset\,Z.$$
\sk
Take any smooth curve $E\subset\PP^2$ of degree $d_E\ges 3$. The genus $g_E$ of $E$ is given by
$$g_E=(d_E-1)(d_E-2)/2.$$
(Note that the vector bundle $E$ in (1.3) is not used in this section.) Set
$$Z_1:=\rho^{-1}(E\mtim\PP^{\1d_2-2})\,\subset\,Z$$
Here $\dim Z_1=\dim Z_2=d_2$. Put
$$X_2:=Z_1\cup Z_2\,\subset\,Z,\q Z'_1:=Z_1\setminus Z_2,$$
with $j_{Z'_1}:Z'_1\into Z_1$ the natural inclusion.
As in (2.2.2), we have the distinguished triangle
$$\DD\Q_{Z_2}\to\DD\Q_{X_2}\to\R(j_{Z'_1})_*\DD\Q_{Z'_1}\buildrel{+1}\over\to.
\leqno(2.3.1)$$
Then
$$\pH^{-d_2}\DD\Q_{X_2}=\DD\Q_{X_2}[-d_2].
\leqno(2.3.2)$$
and we have the long exact sequence of mixed $\Q$-Hodge structures (see \cite{De2}, \cite{Sa2}):
$$\to H_k(Z_2)\to H_k(X_2)\to H^{\rm BM}_k(Z'_1)\to H_{k-1}(Z_2)\to,
\leqno(2.3.3)$$
where $H^{\rm BM}_{\ssb}$ denotes the Borel-Moore homology.
\sk
By the Thom-Gysin sequence (1.3.2) for $\F^{\ssb}=\DD\Q_{B_0}$ with $B_0:=E\mtim\PP^{\1d_2-2}$, we have the long exact sequence of mixed $\Q$-Hodge structures (see Remark after (1.3)):
$$H_k(B_0)\buildrel{c'}\over\rightarrow H_{k-2}(B_0)(1)\to H^{\rm BM}_k(Z'_1)\to H_{k-1}(B_0)\buildrel{c'}\over\rightarrow H_{k-3}(B_0)(1),
\leqno(2.3.4)$$
where $c'$ is the first Chern class of the line bundle $\rho':Z'\to B$.
\sk
Using (2.3.3--4), we can show the isomorphisms of $W$-graded mixed Hodge structures of {\it odd\1} weights
$$\Gr^W_{\rm odd}\1H_k(X_2)=\Gr^W_{\rm odd}\1H^{\rm BM}_k(Z'_1)=\begin{cases}H_1(E)(d_2-1)&\h{if}\,\,\,k=2d_2-1\,,\\ H_1(E)&\h{if}\,\,\,k=2\,,\\ \,0&\h{otherwise,}\end{cases}
\leqno(2.3.5)$$
where
$$\Gr^W_{\rm odd}\1H_k(X_2):=\mopl_{i\in 2\Z+1}\,\Gr^W_iH_k(X_2),\,\,\h{etc.}$$
\skn
Indeed, the first isomorphism of (2.3.5) follows from (2.3.3). For the second, we have
$$\Gr^W_{\rm odd}\1H_k(B_0)=\begin{cases}H_1(E)(j)&\h{if}\,\,\,k=2j+1,\,\,j\in[0,d_2-2],\\ \,0&\h{otherwise,}\end{cases}$$
and (2.3.4) implies the isomorphisms
$$\aligned\Gr^W_{\rm odd}\1H^{\rm BM}_{2d_2-1}(Z'_1)&=\Gr^W_{\rm odd}\1H_{2d_2-3}(B_0)(1)=H_1(E)(d_2-1),\\ \Gr^W_{\rm odd}\1H^{\rm BM}_2(Z'_1)&=\Gr^W_{\rm odd}\1H_1(B_0)=H_1(E),\endaligned$$
where the other $\Gr^W_{\rm odd}\1H^{\rm BM}_k(Z'_1)$ vanish. So (2.3.5) follows.
\sk
Set
$$X:=X_1\times X_2,\q d:=d_1+d_2\,(=\dim X),$$
where $X_1$, $d_1$ are as in (2.1). Then
$$\pH^{-j}\DD\Q_X=\begin{cases}\DD\Q_{X_1}[-d_1]\boxtimes\DD\Q_{X_2}[-d_2]&\h{if}\,\,\,j=d\,,\\ \,0&\h{otherwise.}\end{cases}$$
\sk
To simplify the argument, we assume the following:
$$d_1=2.$$
By the last assumption in (2.1), we have
$$H_1(X_1)=H_3(X_1)=0.$$
The following morphisms are surjective (using (2.1.2)):
$$c_1(D),c_1(D'):H_2(X_1)\onto H_0(X_1).$$
So there is only a {\it difference in the rank\1} for
$$c_1(D),c_1(D'):H_4(X_1)\to H_2(X_1).
\leqno(2.3.6)$$
\sk
We have very ample line bundles $\Lc,\Lc'$ on $X$ defined by
$$\Lc:={\rm pr}_1^*\Lc_1\otimes{\rm pr}_2^*\Lc_2,\q \Lc':={\rm pr}_1^*\Lc'_1\otimes{\rm pr}_2^*\Lc_2.$$
Here ${\rm pr}_i$ denotes the $i\1$th projection from $X_1\times X_2$, $\Lc_1$ is the very ample line bundle on $X_1$ corresponding to $D$ (similarly for $\Lc'_1$ with $D$ replaced by $D'$, assuming $D'-D$ effective and base-point-free), and $\Lc_2$ is a very ample line bundle on $X_2$.
\sk
We now show that condition~(5) holds if $j=j_0$, $k=k_0$ with
$$j_0:=d+1,\q k_0:=(2-d_1)+(d_2-2)=d_2-d_1\,\,(=d_2-2),
\leqno(2.3.7)$$
where $k_0\ges 2$ if $d_2\ges d_1+2\,\,(=4)$. The number $d_2-2$ appears here, since we have in (2.3.5)
$$H_2(X_2)=H^{-2}(X_2,\DD\Q_{X_2})=H^{d_2-2}(X_2,\pH^{-d_2}\DD\Q_{X_2}).$$
We will show that the {\it even} weight part can be neglected effectively in condition~(5) if $g_E\gg 0$. In the notation of (4), set
$$\Gr^W_{\rm even}\1H^k_{(d)}(X)^{\Lc}:=\mopl_{i\in2\Z}\,\Gr^W_iH^k_{(d)}(X)^{\Lc}\q\q(k\in\Z),$$
and similarly for $\Gr^W_{\rm odd}\1H^k_{(d)}(X)^{\Lc}$, $\Gr^W_{\rm even}\1H^k_{(d)}(X)_{\Lc}$, $\Gr^W_{\rm odd}\1H^k_{(d)}(X)_{\Lc}$. Here $W$ is the weight filtration of the canonical mixed Hodge structure on
$$H^k_{(d)}(X):=H^k(X,\pH^{-d}\DD\Q_X)=H_{d-k}(X).$$
\vskip-1mm\nin
Put
\vskip-5mm
$$\aligned&\q\q\q\q\q\mu^k_{\rm odd}(X,\Lc):=\mu^k_{\rm odd}(X)_{\Lc}+\mu^{k-1}_{\rm odd}(X)^{\Lc}\q\q\q\h{with}\\&\mu^k_{\rm odd}(X)_{\Lc}:=\dim\Gr^W_{\rm odd}\1H^k_{(d)}(X)_{\Lc},\q\mu^{k-1}_{\rm odd}(X)^{\Lc}:=\dim\Gr^W_{\rm odd}\1H^{k-1}_{(d)}(X)^{\Lc}\raise12pt\h{}\raise-5pt\h{},\endaligned$$
and similarly for $\mu^k_{\rm even}(X,\Lc)$, $\mu^k_{\rm even}(X)_{\Lc}$, $\mu^{k-1}_{\rm even}(X)^{\Lc}$. For $\Lc$, $\Lc'$, $k_0$ as above, we then get
$$\delta^{k_0}_{\rm odd}:=\bl|\mu^{k_0}_{\rm odd}(X,\Lc)-\mu^{k_0}_{\rm odd}(X,\Lc')\br|>\delta^{k_0}_{\rm even}:=\bl|\mu^{k_0}_{\rm even}(X,\Lc)-\mu^{k_0}_{\rm even}(X,\Lc')\br|,
\leqno(2.3.8)$$
if $g_E\gg 0$. Indeed, $\delta^{k_0}_{\rm odd}$ is strictly positive by (2.1.4), (2.3.5), and is proportional to $g_E$ by Lemma below via the inclusion (using the K\"unneth formula):
$$H^{2-d_1}(X_1,\pH^{-d_1}\DD\Q_{X_1})\otimes H^{d_2-2}(X_2,\pH^{-d_2}\DD\Q_{X_2})\into H^{k_0}(X,\pH^{-d}\DD\Q_X),
\leqno(2.3.9)$$
when the curve $E$ is changed. On the other hand, $\delta^{k_0}_{\rm even}$ in (2.3.8) is independent of $g_E$. So condition~(5) holds if $g_E\gg 0$, see also Remarks~(i--iii) below for more precise arguments. This finishes the proof of Theorem~1.
\msn
{\bf Lemma.} {\it In the above notation, $\delta^{k_0}_{\rm odd}$ in $(2.3.8)$ is proportional to $g_E$, and $\delta^{k_0}_{\rm even}$ remains invariant under the change of the plane curve $E\subset\PP^2$.}
\msn
{\it Proof.} The {\it odd\1} weight part of $H_{\ssb}(X_1)$ is
$$\Gr^W_{-1}H_2(X_1)=H_1(D),$$
and the actions of $c_1(D),c_1(D')$ on it vanish. We get an {\it even} weight part of $H_{\ssb}(X)$ via the K\"unneth formula by taking the tensor product of this {\it odd\1} weight part with the {\it odd\1} weight part of $H_{\ssb}(X_2)$ which is isomorphic by (2.3.5) to
$$H_1(E)\oplus H_1(E)(d_2-1).
\leqno(2.3.10)$$
The action of $c_1(\Lc_2)$ on the latter odd weight part vanishes. Hence the actions of $c_1(\Lc),c_1(\Lc')$ vanish on the above tensor product, which is called the {\it odd-odd\,} weight part. (Here we consider the actions on the graded pieces of the weight filtration $W$. Note that any morphism of mixed Hodge structures is {\it strictly compatible} with the weight filtration, and the kernel and cokernel commute with the passage to the {\it graded quotients} of the weight filtration, see \cite{De1}.) So the contribution of this {\it odd-odd\,} weight part vanishes by taking the difference between $\mu^{k_0}_{\rm even}(X,\Lc)$ and $\mu^{k_0}_{\rm even}(X,\Lc')$. This shows the invariance of $\delta^{k_0}_{\rm even}$ under the change of the curve $E$, since the {\it even-even} weight part is clearly independent of $E$.
\sk
As for $\delta^{k_0}_{\rm odd}$ in (2.3.8), we see that the contribution of the tensor product of the even weight part of $H_{\ssb}(X_1)$ with (2.3.10), that is, the {\it even-odd\,} weight part, is proportional to $g_E$ (where the action of $c_1(\Lc_2)$ vanish on (2.3.10)). Note that only $H_1(E)$ in (2.3.10) contributes here, and there is no contribution of $H_1(E)(d_2{-}1)$ for a reason of degree (since $k_0=d_2-d_1$ in (2.3.7)), see also (2.3.5--6). So it remains to consider the {\it odd-even} weight part. We see that the odd weight part of $H_{\ssb}(X_1)$ does not contribute to $\delta^{k_0}_{\rm odd}$ in (2.3.8) by taking the tensor product with the even weight part of $H_{\ssb}(X_2)$, since there is a {\it difference of actions} only on the {\it even} weight part of $H_{\ssb}(X_1)$ as is explained above, see also (2.3.6). This finishes the proof of Lemma.
\msn
{\bf Remarks.}~(i) Some part of the above Lemma can be avoided if we assume, for instance, $Y=\PP^1\times\PP^1$ with $D\subset\PP^1\times\PP^1$ the diagonal so that $H^1(D)=0$ as in Remark after (2.1).
\ms
(ii) Setting $\delta_2:=d_2-4\ges 0$, we have
$$\aligned\dim H^k(Z_2,\pH^{-d_2}\DD\Q_{Z_2})&=\begin{cases}6-|k|+\delta_2&\h{if}\,\,\,\,|k|-\delta_2=2\,\,\,\,\h{or}\,\,\,\,4,\\6&\h{if}\,\,\,\,|k|\les\delta_2,\,k-\delta_2\in 2\1\Z,\\0&\h{otherwise,}\end{cases}\\
\dim H^k(Z'_1,\pH^{-d_2}\DD\Q_{Z'_1})&=\begin{cases}2g_E&\h{if}\,\,\,\,k=-\delta_2{-}3\,\,\,\,\h{or}\,\,\,\,\delta_2{+}2,\\1&\h{if}\,\,\,\,k=-\delta_2{-}3\pm 1\,\,\,\,\h{or}\,\,\,\,\delta_2{+}2\pm 1,\\0&\h{otherwise.}\end{cases}\endaligned$$
Moreover the following morphisms appearing essentially in (2.3.3) are injective:
$$H^k(Z'_1,\pH^{-d_2}\DD\Q_{Z'_1})\to H^{k+1}(Z_2,\pH^{-d_2}\DD\Q_{Z_2})\q\q(k=d_2{-}2\pm 1).$$
(Note that $d_2{-}2\pm 1=\delta_2{+}2\pm 1$.) This can be proved by using a long exact sequence like (2.3.3) with $Z_2$, $X_2$ respectively replaced by $Z_1\setminus Z'_1$, $Z_1$, together with a morphism from this sequence to (2.3.3). (Here we also study a similar sequence for the $\PP^1$-bundle over $\PP^{\1d_2-2}$.) We then get
$$H^{d_2-1}(X_2,\pH^{-d_2}\DD\Q_{X_2})=0.$$
So the contribution to $\delta^{k_0}_{\rm even}$ in (2.3.8) via the following inclusion vanishes:
$$H^{-d_1}(X_1,\pH^{-d_1}\DD\Q_{X_1})\otimes H^{d_2-1}(X_2,\pH^{-d_2}\DD\Q_{X_2})\into H^{k_0-1}(X,\pH^{-d}\DD\Q_X),$$
where $\mu^{k_0-1}_{\rm even}(X)^{\Lc}-\mu^{k_0-1}_{\rm even}(X)^{\Lc'}$ is involved. \sk
We have also the vanishing of the contribution to $\delta^{k_0}_{\rm even}$ via the inclusion (2.3.9), where $\mu^{k_0}_{\rm even}(X)_{\Lc}-\mu^{k_0}_{\rm even}(X)_{\Lc'}$ is involved. This follows from Remark~(iii) below together with the hard Lefschetz property of the action of $c_1(\Lc_2)$ on $H^{\ssb}(Z_2)$ and the commutativity of the morphisms in (2.3.3) with the action of $c_1(\Lc_2)$. (The argument is rather delicate in the case $d_2=4$, where we need also the primitive decomposition together with (2.3.6).)
\sk
For (2.3.8) it is then sufficient to assume $d_E\ges 3$ so that $g_E\ges 1$, since $\dim H^1(E)=2\1g_E$.
\ms
(iii) Let $A_{\ssb}$, $B_{\ssb}$ be graded vector spaces having the actions $c':A_i\to A_{i+1}$, $c'':B_i\to B_{i+1}$ ($i\in\Z$), and satisfying the following conditions for some integers $p,q$ :
$$A_i=0\,\,(\forall\,i\notin[p,p+2]),\q c''(B_i)=B_{i+1}\,\,(\forall\,i\in[q,q+2]).$$
Set $c:=c'\otimes{\rm id}+{\rm id}\otimes c''$ on $C_{\ssb}:=A_{\ssb}\otimes B_{\ssb}$.
Then $c(C_{p+q+2})=C_{p+q+3}$.

\end{document}